\nonstopmode
\documentclass[a4paper,10pt]{amsart}
\usepackage{amsmath,mathrsfs,amssymb,amsthm,amscd}
\usepackage{a4wide}
\usepackage[]{fontenc}
\usepackage[all]{xy} 
\usepackage{graphicx}
\usepackage{color}

\usepackage[T1]{fontenc}
\usepackage[utf8]{inputenc}
\usepackage{amsfonts}
\usepackage{hyperref}
\usepackage{listings}

\newcommand{\belyi}{{\boldsymbol \beta}}
\newcommand{\BD}{\textbf{D}}

\newcommand{\BS}{\textbf{S}}

\newcommand{\CC}{\mathcal{C}}

\newcommand{\C}{\mathcal{C}}
\newcommand{\E}{\mathcal{E}}
\renewcommand{\O}{\mathcal{O}}
\newcommand{\NZ}{\mathbb{N}}
\newcommand{\GZ}{\mathbb{Z}}
\newcommand{\CZ}{\mathbb{C}}
\newcommand{\RZ}{\mathbb{R}}
\newcommand{\QZ}{\mathbb{Q}}
\newcommand{\OH}{\mathbb{H}}
\newcommand{\Pro}{\mathbb{P}^1}
\newcommand{\G}{\Gamma}
\newcommand{\g}{\gamma}
\newcommand{\PSL}{PSL}
\newcommand{\Div}{\operatorname{Div}}
\newcommand{\hei}{\operatorname{ht}}

\newtheorem{satz}{Theorem}[section]

\newtheorem{prop}[satz]{Proposition}
\newtheorem{lemma}[satz]{Lemma}
\newtheorem{defin}[satz]{Definition}
\newtheorem{bem}[satz]{Remark}

\newenvironment{beweis}{\par\pagebreak[2]\noindent{\it Proof: }}{
 \hfill $\Box$ \medskip}
\newcommand{\mat}[4]{{\left(\begin{smallmatrix}#1&#2 \\ #3&#4\end{smallmatrix}\right)}}

\numberwithin{equation}{section}

\renewcommand{\Re}{\operatorname{Re}}
\renewcommand{\Im}{\operatorname{Im}}

\begin{document}
 
\markright{Busch, K\"uhn, Posingies, \today}
 \title{On scattering constants for a   non-con\-gruence subgroup }
\author{Vincenz Busch, Ulf K\"uhn, Anna Posingies} 

\date{\today} \email{vinzb@gmx.net \\kuehn@math.uni-hamburg.de \\ anna.posingies@math.uni-hamburg.de}
\address{Department Mathematik (AZ)\\Universit\"at Hamburg\\
  Bundesstrasse 55\\D-20146 Hamburg}
%\classification{14G40,11G05,11g30}
\keywords{Arakelov theory, arithmetic intersection numbers, modular curves, elliptic curves, non-congruence subgroups, Eisenstein series}
\thanks{{\bf MSC 2000:} 14G40, 11g05, 11g50, 11M36} 
%Arithmetic varieties and schemes; Arakelov theory; heights, 
%Elliptic curves, 
%Structure of modular groups and generalizations; arithmetic groups, 
%Selberg zeta functions and regularized determinants; applications to spectral theory, 
%Dirichlet series, Eisenstein series, etc. Explicit formulas

\maketitle
\begin{abstract} 
  Scattering constants are special values of Dirichlet series
  associated to non-holomorphic Eisenstein series.  In this paper we
  give closed formulas for the scattering constants related to a
  non-congruence subgroup obtained via a Belyi map of an elliptic
  curve.
\end{abstract}

 \section*{Introduction} 
 Let $\Gamma \subset SL_2(\GZ)$ be  a finite index subgroup. To every cusp $S_j$ of $\Gamma$ we can associate a non-holomorphic Eisenstein series $E_j^{\G}(z,s)$ and for two cusps $S_j$ and $S_k$ the scattering constant $C^{\G}_{jk}$ is defined by the constant term of a Dirichlet series coming from the Fourier expansion of $E_j^{\G}(z,s)$ in the cusp $S_k$. \\
 Given an algebraic curve $C$ and a Belyi map $\belyi$, i.e. a
 morphism $\belyi: C \rightarrow \Pro$ that is ramified in at most
 three points, for simplicity we assume that both are defined over
 $\QZ$, then there exists a subgroup $\G \subset \G(2)$ such that
 $C(\CZ)\cong \overline{\G \setminus \OH}$.  The cusps of $\G$
 correspond to the ramification points of $(C, \belyi)$. Arakelov
 theory gives us an expression for the N\'eron Tate height (see
 theorem \ref{satzKuehn} below for more details)
$$\hei_{NT}(\text{ ramification points }) = \text{weighted sum of } C_{jk}^\Gamma\, 's+ \text{algebraic term}.$$
Our main result concerns the elliptic curve  $E: y^2 = x^3+5x+10$ and
the non-congruence subgroup $\G_E$ for the
Belyi map
%	\begin{align}
$	\belyi_E: E  \to  \Pro$   
%	\label{eq:belyipair}
%	\end{align}	
     given by $\belyi_E(x,y)= \frac{y(x-5)+16}{32}.$ This curve is
        referred to as 400H1 in Cremona's tables \cite{cremona}, its Mordell Weil
        group has rank one and is generated by $P=(1,4)$.  The zero
        element $\mathcal O$ and $P$ as well as $-P$ correspond to
        cusps of the group $\G_E$.
        \medskip \\
        {\bf Theorem:} {\it For the non-congruence subgroup $\G_E$ as
          above the scattering constants $C^{\G_E}_{\mathcal
            O,\mathcal O}, C^{\G_E}_{P,P}$ and $C^{\G_E}_{P,-P}$ are
          given by
	\begin{align*}
	 C^{\G_E}_{\mathcal O,\mathcal O} &=  \frac{1}{30}\left(C^{\G(1)}-\frac{1}{\pi}\left( 14 \log(2)+6\log(5) \right)\right) \\
	 C^{\G_E}_{P,P} &=  \frac{1}{120}\left(4 C^{\G(1)}+ \frac{1}{\pi}\left(-131\log(2) +15\log(5)-60 \hei_{NT}(P) \right)\right)\\
	 C^{\G_E}_{P,-P} &=   \frac{1}{120}\left(4 C^{\G(1)}+\frac{1}{\pi}\left( -71\log(2)+ 60 \hei_{NT}(P) \right)\right),
	\end{align*}
	where 
	\begin{align}
	  C^{\Gamma(1)} = -\frac{6}{\pi}\left( 12\zeta'(-1)-1+\log(4\pi)\right) \label{eq:SKG1}
	\end{align}
	is the unique scattering constant for the full modular group $\G(1)$.\\
	All other scattering constants for $\G_E$ are $\QZ$-linear combinations of the four scattering constants given above. The coefficients depend on the ramification data coming from $\belyi_E$.} 

      {\bf Remark:} {\it Previous to these explicit expressions, nearly
        no formulas for scattering constants have been available.  For
        Eisenstein series coming from certain congruence subgroups of
        $SL_2(\GZ)$ formulas for the scattering matrices are known
        (\cite{hejhal}, \cite{huxley}). Beside that, Venkov studied
        cycloidal groups \cite{venkov}. From these results scattering
        constants can be deduced.}

	\section{Eisenstein Series}
	\label{kap:ESR}

	Let $\OH=\left\lbrace z\in \CZ \,|\, \Im(z)>0 \right\rbrace $
        be the upper half plane and $\Gamma(1)=\PSL_2(\GZ)$. Then
        $\Gamma(1)$ acts on $\OH$ by the M\"obius
        transformation. These action can be extended to $
        \overline{\OH}=\OH \cup \QZ \cup \infty$. Let $\Gamma \subset
        \Gamma(1)$ be a finite index subgroup then $\overline{\QZ}=\QZ
        \cup \infty$ is divided into finitely many equivalence classes
        with respect to the action of $\Gamma$; the classes are called
        cusps of $\Gamma$. We will use the word cusp as well for a
        representative of a cusp. Let $S_j\in \overline{\QZ}$ be a
        cusp and $\Gamma_j$ its stabilizer in $\G$. For $S_j$ there is
        a $\gamma_j \in \Gamma(1)$ with $\gamma_j(\infty)=S_j$ and a
        $b_j \in \NZ$, such that
	\[  \sigma_j^{-1} \Gamma_j \sigma_j = \left\langle \mat{1}{1}{0}{1} \right\rangle \quad \text{with } \sigma_j=\gamma_j \cdot \mat{\sqrt{b_j}}{0}{0}{1/\sqrt{b_j}}.  \] 
	The number $b_j$ is called the width of the cusp $S_j$.

	\begin{defin} Let $\G \subset \G(1)$ be a finite index subgroup. For each cusp $S_j$ there is a
	\emph{non-holomorphic Eisenstein series $E^{\G}_j(z,s)$}, which for
	$z \in \OH$, $s \in \CZ$ and $\Re s >1$ is defined by the convergent series
	$$
	E^{\G}_j(z,s) = \sum_{\sigma \in \Gamma_j \setminus \Gamma} \Im
	\left( \sigma_j^{-1} \sigma(z)\right)^{s}= b_j^{-s}\sum_{\sigma \in \Gamma_j \setminus \Gamma} \Im
	\left( \gamma_j^{-1} \sigma(z)\right)^{s}.
	$$
	\end{defin}

	\noindent {\bf Properties:}  Let us recall some facts on the theory of
	Eisenstein series; the standard reference is \cite{kubota}. The function $E^{\G}_j(z,s)$ has a meromorphic
	continuation to the $s$-plane, with a simple pole in $s=1$ with
	residue $3/( \pi \cdot [ \Gamma(1): \Gamma])$.  For all $\gamma\in
	\Gamma$ we have $E^{\G}_j(\gamma(z),s)=E^{\G}_j(z,s)$. 
	The Fourier expansion of $E^{\G}_j(z,s)$ at the cusp $S_k$ is given
	by
	$$ 
	E^{\G}_j(\sigma_k(z),s)= \delta_{jk}\cdot y^s+\pi^{1/2}\frac{\Gamma\left( s-\frac{1}{2}\right) }{\Gamma(s)} \cdot \phi^{\G}_{jk}(s)\cdot y^{s-1} + \sum_{m \neq 0}a_m(y,s)e^{2\pi i m x}
	$$
	where $z=x+iy$ and $\Gamma(s)$ is the Gamma function. Furthermore we have
	\begin{align}
	 \phi^{\G}_{jk}(s)= \frac{1}{(b_jb_k)^s}\sum_{c>0}r_{jk}^{\Gamma}(c)\frac{1}{c^{2s}}  \label{eq:kleinphi}
	\end{align}
	and
	\begin{equation}
	r_{jk}^{\Gamma}(c) =\# \left\lbrace d \mod b_kc \,|\, \exists \, \mat{*}{*}{c}{d} \in \gamma_j^{-1} \Gamma \gamma_k  \right\rbrace  .
	\label{eq:defrjk}
	\end{equation}
	Then the \emph{scattering matrix } 
	$$
	\Phi_\Gamma(s) = \left(\pi^{1/2}
	\frac{\Gamma(s-1/2)}{\Gamma(s)} \cdot  \phi^{\G}_{jk}(s) \right)_{j,k}
	$$
	is symmetric. Note that all the
	coefficients of the scattering matrix are Dirichlet series in a general
	sense. They have a meromorphic continuation with a simple pole in $s=1$
	of residue $3/( \pi \cdot [ \Gamma(1): \Gamma])$.

	\begin{defin} \label{def:Cjk} For all pairs $j,k$ we define
	the \emph{scattering constant} $C^{\G}_{jk}$ to be the constant term at
	$1$ of the Dirichlet series $(\Phi_\Gamma)_{jk}(s)$, i.e.
	\begin{align}
	C^{\G}_{jk}:= \lim_{s \to 1} \left(\Phi_\Gamma(s)_{j,k} -
	\frac{3/(\pi \cdot [\Gamma(1):\Gamma])}{s-1}   \right).
	\end{align}
	\end{defin}

%	The scattering matrices are known for congruence subgroups (\cite{hejhal} and \cite{huxley}).
	We will need the values of the scattering constants for $\Gamma(1)$ and $\G(2)$.
	The group $\G(1)$ has one cusp, hence one scattering constant and it has already been introduced in \eqref{eq:SKG1}.
	For the group $\G(2)$ there exist only two different scattering constants, although $\G(2)$ has three cusps. These constants are \begin{align}
		 C^{\G(2)}_{a} & = -\frac{1}{3\pi} \left( 36\zeta'(-1)-3+3\log(4\pi)+7\log(2) \right) \notag \\
			& = \frac{1}{6}C^{\Gamma(1)}-\frac{7}{3\pi}\log(2)  \label{eq:SKG2gleich}\\
		 C^{\G(2)}_{b} & = -\frac{1}{3\pi} \left( 36\zeta'(-1)-3+3\log(4\pi)+\log(2) \right) \notag \\
			& =  \frac{1}{6}C^{\Gamma(1)}-\frac{1}{3\pi}\log(2), \label{eq:SKG2ungleich}
		\end{align}
		where the first case is the one with $S_j=S_k$ and in the second case we have $S_j \neq S_k$. For a calculation see e.g. \cite{posingies}. \medskip
%		where the first case is the one with $S_j=S_k$ and in the second one we have $S_j \neq S_k$. The formulas for the scattering constants can be extracted from \cite{hejhal} or \cite{huxley}. Furthermore, they are not difficult to calculate directly from the definition as it is done in \cite{posingies}.
	
		The scattering constants for a group can be constructed from the constants for a subgroup. Conversely, the knowledge of scattering constants for a group gives us some information about sums of scattering constants for subgroups. Take the group $\G(2)$, then we get:
	\begin{prop}
		\label{pr:sumsc}
	Let $\G$ be a finite index subgroup of $\G(2)$, $S_j^{\G(2)}$ and $S^{\G(2)}_{k'}$ two cusps of $\G(2)$. Then 
	 \begin{align}
	\sum_{S_i^{\Gamma}\subset S_j^{\Gamma(2)}}\frac{b_i}{2} C_{ik}^{\Gamma}=C_{jk'}^{\Gamma(2)}- 
	\frac{1}{2\pi [\Gamma(2):\Gamma] } \sum_{S_i^{\Gamma}\subset S_j^{\Gamma(2)}}\frac{b_i}{2}
	\log\left( \frac{b_ib_k}{4}\right),
		\label{eq:sumsc}
	\end{align}
	where we sum over a system of representatives $\lbrace S_i^{\G} \rbrace$ of cusps of $\G$ such that $S_i^{\Gamma}\sim_{_{\G(2)}} S_j^{\G(2)}$ and $S_k^{\G}$ is any cusp of $\G$ with $S_k^{\G}\sim_{_{\G(2)}} S_{k'}^{\G(2)}$. The $b_i$ and $b_k$ denote the widths of the cusps. 
	\end{prop}
		\begin{beweis}
		The formula is a consequence of the relation of Eisenstein series 
		\begin{align*}
		 2^{s}E_j^{\G(2)}(z,s) =\sum_{\g \in \G_j \setminus \G(2)}\Im (\g_j \g (z))^s &=
		\sum_{-\frac{d}{c}\in S_j^{\G(2)}}\frac{\Im(z)^s}{|cz-d|^{2s}} \\
		&= \sum_{S_i^{\G} \subset S_j^{\G(2)}} \sum_{-\frac{d}{c}\in S_i^{\G}}\frac{\Im(z)^s}{|cz-d|^{2s}} \\
		&= \sum_{S_i^{\G} \subset S_j^{\G(2)}} b_i^s E_i^{\G}(z,s).
		\end{align*}
		From the implied relation for the constant terms of the Fourier expansions we conclude the claim by a straight forward calculation.

%		A relation for every coefficient in the Dirichlet series \eqref{eq:kleinphi} follows, e.g.
%		$$\frac{1}{2}r_{jk'}^{\G(2)}(c)= \frac{1}{b_k}\sum_{S_i^{\G} \subset S_j^{\G(2)}} r_{ik}^{\G}(c) \quad \forall c\in \NZ. $$
%		With this last identity the formula from the statement can be calculated. 
		\end{beweis}

\section{Belyi's Theorem,  N\'eron Tate Heights and Arakelov Theory}

	In 1979 G. Belyi proved in \cite{belyi}
	\begin{satz}
	Let $C$ be a non-singular algebraic curve defined over a number field. Then there exists a finite morphism
	$\belyi: C \longrightarrow \Pro$ with at most the three critical values $0,1, \infty$.
	\end{satz}

	This was the missing step for the following equivalences.
	\begin{satz}
	\label{sz:equivbelyi}
	 Let $C$ be a non-singular algebraic curve over $\CZ$. Then the following are equivalent
		\begin{enumerate}
		 \item  The curve $C$ is defined over a number field.
		 \item  There exists a finite morphism $\belyi_C: C \longrightarrow \Pro$ with at most the three critical values $0,1, \infty$.
		 \item  There is a subgroup $\Gamma_C \subset \Gamma(2)$ such that $ C(\CZ) \cong \overline{\G_C \setminus \OH} .$ 
		\end{enumerate}
	\end{satz}
		\begin{beweis}
		See e.g. \cite{birch}, \cite{bost} or \cite{serre}.
		\end{beweis}

	\begin{defin}
          A pair consisting of a curve and a map with the properties
          from $(ii)$ is called Belyi pair. The map alone is a Belyi
          map.
	\end{defin}
	Some data from the Belyi pair have a direct counterpart in the
        group $\G_C$. There is a 1-1 correspondence between the cusps
        of $\G_C$ and the ramification points of the Belyi pair. The
        widths of the cusps resemble the ramification orders.
	\begin{defin}
          We will call the ramification points of a Belyi pair
          cusps. A divisor on $C$ which has only support in the cusps
          of $(C, \belyi_C)$ is called cuspidal divisor.
	\end{defin}

	Belyi pairs allow to formulate the following theorem
        concerning N\'eron Tate heights.
	\begin{satz}\label{satzKuehn}
          Let $\belyi: K \rightarrow \Pro$ be a Belyi map for an
          algebraic curve $C$ over $\QZ$ with induced Belyi
          uniformization $C(\CZ)\equiv \overline{\Gamma_C \setminus
            \OH}$. Let $D = \sum_j n_jS_j$ and $D' = \sum_k m_kS_k$ be
          two cuspidal divisors of degree $0$.
%, i.e., $D, D' \in \Div_0(C)$.  
Then the N\'eron Tate height pairing of $D$ and
          $D'$ is given by:
	\begin{equation}
	\left\langle D,D' \right\rangle_{NT} = - \sum_{p \text{ prim}}\delta_p \log(p)- 2\pi \sum_{j,k} n_j m_k C^{\Gamma_C}_{jk}.
	\label{eq:kuehn}
	\end{equation}
	The coefficients $\delta_p$ are rational numbers, that are
        explicitly computable (see \eqref{EQfin}) and
        $C^{\Gamma_C}_{jk}$ is the scattering constant for the cusps
        $S_j$ and $S_k$ from $\G_C$.
	\label{sz:kuehn}
	\end{satz}
	\begin{beweis}
	See \cite{kuehn}.
	\end{beweis}
	
	The above coefficients $\delta_p$, as we will explain now,
        are given via  local intersection numbers.

	Let $C/\QZ$ be an algebraic curve and $\C/\GZ$ be a 
        proper regular model.  Let $\BD_1, \BD_2 \in \Div(\C)$ be two
        prime divisors with no common components and $x \in \BD_1\cap
        \BD_2$.  Fix local parameters $f_1, f_2$ for $\BD_1$ and
        $\BD_2$.  We define the local intersection number at $x$ to be
	\[i_x(\BD_1,\BD_2) = \ell_{\O_{\C,x}}(\O_{\C,x})/(f_1,f_2).\]
       	Further, we define the total local intersection number of  $\BD_1$ and $\BD_2$ at a
        prime $p$ to be
	\[i_p(\BD_1,\BD_2) = \sum_{x \in \BD_1\cap \BD_2 \cap \C_p} i_x(\BD_1,\BD_2) [k_x:k_p]. \]
	For two divisors $\BD = \sum n_j \BD_j$, $\BD' = \sum m_k
        \BD'_k$ with no common components we define by linearity
	\[i_p(\BD,\BD') = \sum_{j,k} n_j m_k i_p(\BD_j,\BD'_k). \]
Now we  are looking at the different fibers of the scheme $\C/\GZ$ at once.
	\begin{defin}
 For two divisors $\BD,\BD' \in \Div(\C)$ with no common
          component we define the intersection number at the finite places by
	\[(\BD,\BD')_{fin} = \sum_{p} i_p(\BD,\BD') \log(p).\]
	\end{defin}

	From now on $\BD$ denotes the Zariski closure of a divisor $D$
        on $C$.  The group of divisors on $C$ with degree zero will be
        denoted by $\Div_0(C)$ and $\Div_p(\C)$ is the set of all
        divisors supported on $\C_p$, here $\C_p=\C \times k_p\ (k_p =
        \GZ/p \GZ)$ denotes the special fiber of $\C$ at the place
        $p$.
	
	\begin{lemma}
          There exists a unique linear map
          \[ \Phi_p: \Div_0(C) \rightarrow \QZ \otimes \Div_p(\C)
          /(\QZ \otimes \C_p),\] such that for all $D \in \Div_0(C)$
          the divisor $\BD + \Phi_p(D)$ is orthogonal to $\Div_p(\C)$.
	\end{lemma}
	\begin{beweis}
	See e.g. \cite{hriljac}.
	\end{beweis}

	The correction divisor $\Phi$ for $D\in \Div_0(C)$ is defined by
	\[ \Phi(D) = \sum_p \Phi_p(D). \]
%	and usually denote $\CD= \BD + \Phi(D)$.
	
\begin{defin}	We define for two divisors $D,D' \in \Div_0(C)$ with no common components
	\begin{equation*}
	(D,D')_{fin} = (\BD + \Phi(D),\BD' + \Phi(D'))_{fin}.\end{equation*}
\end{defin}

	From the definitions we see
\begin{align}\label{EQfin}
	(D,D')_{fin} = \sum_p \delta_p \log(p). 
\end{align}
	These $\delta_p$ are exactly the $\delta_p$ of theorem \ref{satzKuehn}, whenever 
$D$ and $D'$ have no common components. Since the N\'eron Tate height pairing 
vanishes on divisors of rational functions, we may reduce the general case to the
the above situation by replacing $D'$ with $D' +\operatorname{div} (f)$ for an suitable rational function $f$
on $C$.  
	
\begin{bem}
  As a consequence of the functoriality of the intersection number at
  the finite places w.r.t. pull-back morphisms (see \cite{lang}), the
  quantity $(D,D')_{fin}$ does not depend on the particular chosen
  regular model.
		\end{bem}
	
\section{An elliptic Belyi pair and its scattering constants}

	Now, we will focus on one particular Belyi pair considered by Elkies in \cite{elkies}.	
	 
	\begin{prop}
	 The elliptic curve 
		\begin{align}
		 E: y^2=x^3+5x+10 \label{eq:gleichE} 
		\end{align}
	together with the map 
		\begin{align}
		\belyi: \quad E & \longrightarrow \Pro \label{eq:belyiE} \\
		(x,y) &\longmapsto \frac{y(x-5)+16}{32} \notag
		\end{align}
	form a Belyi pair.
	\end{prop}
		\begin{beweis}
                  Regard the map $(x,y) \mapsto y(x-5)$. It is easy to
                  check, that the critical values are $\lbrace \infty
                  ,\pm 16 \rbrace$. The map \eqref{eq:belyiE} above
                  normalizes the critical values to $0$, $1$ and
                  $\infty$.
		\end{beweis}

	\begin{prop}
          The Mordell Weil group $E(\QZ)$ of the elliptic curve
          defined via \eqref{eq:gleichE} is of rank one with trivial
          torsion. The point $S_1=(1,4)$ generates $E(\QZ)$.  The
          group $\G_E$ associated with the Belyi pair $(E, \belyi_E)$
          is a non-congruence subgroup.
	\end{prop}
		\begin{beweis}
                  For the rank, torsion, and the generator look at
                  Cremona's tables in \cite{cremona}. The curve is
                  listed under 400H1.  The point $S_1$ generates an
                  infinite group.  Hence, the cuspidal divisor
                  $S_1-\mathcal O$ is not torsion in the Picard
                  group. The theorem of Manin and Drinfeld
                  \cite{drinfeld}, \cite{manin} (see also
                  \cite{elkik}) then shows that $\G_E$ is
                  non-congruence.
		\end{beweis}	

\begin{bem}
It was shown in \cite{posingies} that $\G_E=\varphi^{-1}(Stab_{S_5}(5))$ with $\varphi: \G(2) \rightarrow S_5$ given by the images of the generators of $\G(2)$: $\varphi(\g_0)=(1235)$ and  $\varphi(\g_1)=(1234)$. But such a description of $\G_E$ is not needed to achieve the results of this paper.
\end{bem}

The ramification points, i.e. the cusps of the Belyi pair, are  
		\begin{align}
                  & S_0=\mathcal O  \notag \\
                  & S_1=(1, 4) , \quad S_2= -4S_1=(6, -16),  \label{eq:rampoints}  \\
                  & S_3=-S_1=(1, -4),\quad \text{ and } \quad
                  S_4=4S_1=(6, 16), \notag
		\end{align}
                where $\mathcal O$ corresponds to the point $\infty$;
                it is the neutral element of $E(\QZ)$.  The point
                $S_0$ lies above $\infty$, $S_1$ and $S_2$ above $0$
                and $S_3$ as well as $S_4$ above $1$.  The notation
                $S_i$ with $i\in \left\lbrace 0,\dots , 4 \right
                \rbrace$ will be used for the divisors of the cusps in
                $\Div(E)$ and also for the corresponding cusps of the
                group $\G_E$.

                The scattering constants $C_{jk}^{\G_E}$ will be
                determined below via an application of theorem
                \ref{satzKuehn} on the curve $E$ defined via
                \eqref{eq:gleichE} together with the additional
                relations coming from proposition \ref{pr:sumsc}.

                In particular, we need certain intersection numbers at
                the finite places, these will be determined in section
                \ref{kap:delta} using an appropriate model.

\section{Calculation of some coefficients $\delta_p$ for $E$}
\label{kap:delta}
	The elliptic curve $E$  defined via \eqref{eq:gleichE} has a minimal proper regular model
over $\GZ$. We denote the Zariski closure of $S_j$ on the model with $\BS_j$.
	
\begin{satz}
 An illustration of the model and the behavior of the divisors associated to the points from \eqref{eq:rampoints} can be seen in figure \ref{AbbModell400H1}.
	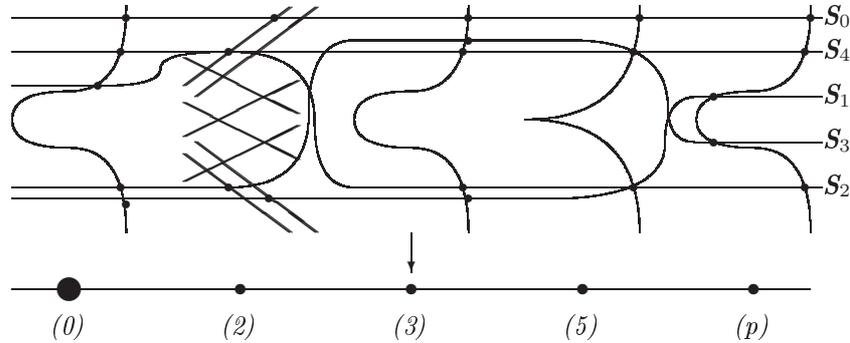
\begin{figure}[ht]
		\centering
\setlength{\unitlength}{.75cm}
\begin{picture}(14.4,6)
%Kurven
%(0)
\qbezier(2,2)(2,3.5)(1,3.5)
\qbezier(1,3.5)(0,3.5)(0,4)
\qbezier(0,4)(0,4.5)(1,4.5)
\qbezier(1,4.5)(2,4.5)(2,6)

%(2)
\put(3,3.4){\line(4,-3){1.86}}
\put(3.2,3.6){\line(4,-3){2.12}}

\put(3,2.9){\line(2,1){2}}
\put(3,4.3){\line(2,-1){2}}
\put(3,3.7){\line(2,1){2}}
\put(3,5.1){\line(2,-1){2}}

\put(3,4.6){\line(4,3){1.86}}
\put(3.2,4.4){\line(4,3){2.12}}

\put(3.02,3.4){\line(4,-3){1.86}}
\put(3.22,3.6){\line(4,-3){2.12}}

\put(3.02,2.9){\line(2,1){2}}
\put(3.02,4.3){\line(2,-1){2}}
\put(3.02,3.7){\line(2,1){2}}
\put(3.02,5.1){\line(2,-1){2}}

\put(3.02,4.6){\line(4,3){1.86}}
\put(3.22,4.4){\line(4,3){2.12}}

\put(3.05,3.4){\line(4,-3){1.86}}
\put(3.25,3.6){\line(4,-3){2.12}}

\put(3.05,2.9){\line(2,1){2}}
\put(3.05,4.3){\line(2,-1){2}}
\put(3.05,3.7){\line(2,1){2}}
\put(3.05,5.1){\line(2,-1){2}}

\put(3.05,4.6){\line(4,3){1.86}}
\put(3.25,4.4){\line(4,3){2.12}}

%(3)
\qbezier(8,2)(8,3.5)(7,3.5)
\qbezier(7,3.5)(6,3.5)(6,4)
\qbezier(6,4)(6,4.5)(7,4.5)
\qbezier(7,4.5)(8,4.5)(8,6)

%(5)
 \qbezier(11,2)(11,4)(9,4)
 \qbezier(11,6)(11,4)(9,4)

%(p)
\qbezier(14,2)(14,3.5)(13,3.5)
\qbezier(13,3.5)(12,3.5)(12,4)
\qbezier(12,4)(12,4.5)(13,4.5)
\qbezier(13,4.5)(14,4.5)(14,6)

%Pfeil
\put(7,2){\vector(0,-1){0.7}}
\put(0,1){\line(1,0){14}}

%Spektrum
\put(1,1){\circle*{0.4}}
\put(.65,.2){\parbox{1cm}{(0)}}

\put(4,1){\circle*{0.2}}
\put(3.65,.2){\parbox{1cm}{(2)}}

\put(7,1){\circle*{0.2}}
\put(6.65,.2){\parbox{1cm}{(3)}}

\put(10,1){\circle*{0.2}}
\put(9.65,.2){\parbox{1cm}{(5)}}

\put(13,1){\circle*{0.2}}
\put(12.65,.2){\parbox{1cm}{(p)}}

%Divisoren

%S0
\put(0,5.8){\line(1,0){14.2}}
\put(14.2,5.7){\parbox[t]{1cm}{$\BS_0$}}

%S1
%\put(6,4.4){\line(1,0){8.2}}
\put(12,4.4){\line(1,0){2.2}}
\qbezier(12,4.4)(11.45,4.4)(11.45,3.6)
%\qbezier(11.45,3.6)(11.45,3.0)(10.9,2.8)
%\qbezier(10.9,2.8)(10.9,2.6)(10,2.4)
\qbezier(11.45,3.6)(11.45,2.6)(9.4,2.6)
\put(0,2.6){\line(1,0){9.4}}

\put(14.2,4.3){\parbox[t]{1cm}{$\BS_1$}}

%S2
\put(6,2.8){\line(1,0){8.2}}
\qbezier(6,2.8)(5.3,2.8)(5.3,4)
\qbezier(5.3,4)(5.3,5.2)(3.8,5.2)
\qbezier(1.6,4.6)(2.6,4.6)(2.6,4.9)
\qbezier(2.6,4.9)(2.6,5.2)(3.8,5.2)
\put(0,4.6){\line(1,0){1.6}}

\put(14.2,2.7){\parbox[t]{1cm}{$\BS_2$}}

%S3
\put(12,3.6){\line(1,0){2.2}}
\qbezier(12,3.6)(11.5,3.6)(11.5,4.2)
\qbezier(11.5,4.2)(11.5,5.4)(9.8,5.4)
\put(6,5.4){\line(1,0){3.8}}
\qbezier(6,5.4)(5.2,5.4)(5.2,4.1)
\qbezier(5.2,4.1)(5.2,2.8)(3.8,2.8)
\put(0,2.8){\line(1,0){3.8}}

%\qbezier(6,3.6)(5.35,3.6)(5.35,3.05)
%\qbezier(5.35,3.05)(5.35,2.5)(4.7,2.5)
%\put(0,2.5){\line(1,0){4,7}}

\put(14.2,3.5){\parbox[t]{1cm}{$\BS_3$}}

%S4
\put(0,5.2){\line(1,0){14.2}}
%\put(7.9,5.2){\line(1,0){6.1}}

\put(14.2,5.1){\parbox[t]{1cm}{$\BS_4$}}

%Schnittpunkte auf p
\put(12.3,4.4){\circle*{0.16}}
\put(12.3,3.6){\circle*{0.16}}
\put(13.9,5.2){\circle*{0.16}}
\put(13.9,2.8){\circle*{0.16}}
\put(14,5.8){\circle*{0.16}}

%Schnittpunkte auf 5
\put(10.9,5.2){\circle*{0.16}}
\put(10.9,2.8){\circle*{0.16}}
\put(11,5.8){\circle*{0.16}}

%Schnittpunkte auf 3
\put(8,2.6){\circle*{0.16}}
\put(8,5.4){\circle*{0.16}}
\put(7.9,5.2){\circle*{0.16}}
\put(7.9,2.8){\circle*{0.16}}
\put(8,5.8){\circle*{0.16}}

%Schnittpunkte auf 2
\put(4.5,2.6){\circle*{0.16}}
\put(3.8,2.8){\circle*{0.16}}
\put(3.8,5.2){\circle*{0.16}}
\put(4.6,5.8){\circle*{0.16}}

%Schnittpunkte auf 0
\put(1.5,4.6){\circle*{0.16}}
\put(2,2.5){\circle*{0.16}}
\put(1.9,5.2){\circle*{0.16}}
\put(1.9,2.8){\circle*{0.16}}
\put(2,5.8){\circle*{0.16}}

\end{picture}

		\caption{The associated divisors on the minimal proper regular model of the curve 400H1}
		\label{AbbModell400H1}
	\end{figure}

	The model fulfills the following properties:
	\begin{enumerate}
	\item 
	All special fibers of $\E$ consist of one irreducible component except the fiber above $(2)$.
	
	\item The fiber above $(2)$ consists of eight irreducible components with the following multiplicities
	\[ \E_2 = \CC_1 + \CC_2 + 2 \CC_3 + 2\CC_4 + 2\CC_5 +2 \CC_6 + \CC_7+ \CC_8\]
	and the intersection matrix of the fiber over $(2)$ is given by table \ref{tab:intmat}. 
	
	\begin{table}[ht] \centering
	\begin{tabular}{c|rrrrrrrr}
	$(\ \cdot\ )$&$\CC_1$ &$\CC_2$ &$2\CC_3$ &$2\CC_4$ &$2\CC_5$ &$2\CC_6$ &$\CC_7$ &$\CC_8$\\
	\hline
	$\CC_1$  &-2  &0   &2   &0   &0   &0   &0   &0 \\
	
	$\CC_2$  &0   &-2  &2   &0   &0   &0   &0   &0 \\
	
	$2\CC_3 $ &2   &2   &-8  &4   &0   &0   &0   &0 \\
	
	$2\CC_4$  &0   &0   &4   &-8  &4   &0   &0   &0 \\
	
	$2\CC_5$  &0   &0   &0   &4   &-8  &4   &0   &0 \\
	
	$2\CC_6$  &0   &0   &0   &0   &4   &-8  &2   &2 \\
	
	$\CC_7 $ &0   &0   &0   &0   &0   &2   &-2  &0 \\
	
	$\CC_8 $ &0   &0   &0   &0   &0   &2   &0   &-2\\ 
	\end{tabular} \medskip
\caption{Intersection matrix for $\E_2$ } \label{tab:intmat}
	\end{table}

	\item There are only three points on the model where two of the divisors $\BS_0, \dots \BS_4$ intersect.
	We denote these points by
	\[ x_{12}=\BS_1 \cap \BS_2\qquad x_{34}=\BS_3 \cap \BS_4\qquad x_{24}=\BS_2 \cap \BS_4.\]
	
	\item
	The divisors $\BS_0, \dots \BS_4$ only intersect with one component of the special fiber above $(2)$ and these intersection numbers are
	 \[ (\BS_0,\CC_1)_{fin} = (\BS_1,\CC_7)_{fin} = (\BS_2,\CC_1)_{fin} = (\BS_3,\CC_8)_{fin} = (\BS_4,\CC_1)_{fin} =1\]
	while all other intersection numbers vanish.
	\end{enumerate}
	\end{satz}
		\begin{beweis} The  Tate algorithm provides an explicit method to obtain 
		the model and its properties $(i),(ii)$.  
Keeping track of the points $S_j$ and their induced divisors $\BS_j$  
in all steps of these calculations gives $(iii)$ and $(iv)$; a detailed calculation
 is in \cite{busch}.  
		\end{beweis}

	Now, we need to calculate two more ingredients to determine the possible $\delta_p$'s for $E$: the local intersection numbers and the correction divisors.

	\medskip	
	\noindent {\bfseries Local intersection numbers:}
	We will denote the local intersection number at $x_{jk}$ with $m_{jk}$.
%	In this paper we will calculate only $m_{QT}$, because the other two intersection numbers are computed analogue.
	Explicit calculations show, see e.g. \cite{busch}, p. 79, that the local ring at the point $x_{24}$ is given by
	\[\O_{\E,x_{24}}=\GZ[x,y]/(y^2-x^3-5x-10)_{\left((2,x,y)\right)}\]
	and the local equations of $\BS_2$ and $\BS_4$ are $f_2 = y-16$ and $f_4 = y+16$.
	
	Thus, the module $\O_{\E,x_{24}} /(f_2,f_4)$ is given by
	
	\[\begin{aligned}
	\O_{\E,x_{24}} /(f_2,f_4) &= (\GZ[x,y]/(y^2-x^3-5x-10)_{\left((2,x,y)\right)})/(y-16,y+16)\\
	&=\GZ[x,y]/(y^2-x^3-5x-10,y-16,y+16)_{\left((2,x,y)\right)}\\
	&=\GZ[x]/(-x^3-5x+246,32)_{\left((2,x)\right)}\\
	&=\GZ/(32)_{\left((2)\right)}.
	\end{aligned}\]
	
	The intersection number is therefore given by
	\[m_{24}=i_{x_{24}}(S_2,S_4) = \ell_{\O_{\E,x_{24}}} \O_{\E,x_{24}}/(f_2,f_4) = \ell_{\O_{\E,x_{24}}} \GZ/(32) = 5.\]
	
	Similar calculations for the other two points give
	\[ m_{12}= m_{34} = 1.\]

	\noindent {\bfseries Correction divisors:}
	We will work with cuspidal divisors that are the difference of
        two cusps and we define $D_{jk}=S_j-S_k$, with $j,k \in
        \lbrace 0, \dots 4 \rbrace$. The Zariski closure of $D_{jk}$
        will be denoted with $\BD_{jk}$.  Next we will calculate the
        divisors $\Phi_{jk} = \Phi(D_{jk})$ with $j,k\in \lbrace 0,
        \dots ,4 \rbrace$.
	
	In our example we have $\Phi_{jk} \in \Div_{(2)}(\E)$, because
        the $\BD_{jk}$ are already orthogonal to all the other fibers.
        Writing $\Phi_{jk}$ as a linear combination of the components
        of the fiber above $(2)$
	\[ \Phi_{jk} = n_1 \CC_1 + n_2 \CC_2 +n_3 \CC_3 + n_4 \CC_4 +
        n_5 \CC_5 + n_6 \CC_6 + n_7 \CC_7 + n_8 \CC_8,\] with $n_i \in
        \QZ$ and solving the equations
	\begin{equation}\label{PhiEquation}
	(\BD_{jk}+\Phi_{jk},\CC_l)_{fin} = 0,\quad 1\leq l \leq 8
	\end{equation}
	we get representatives for the coefficients $n_i$'s; observe $
	\Phi_{jk}$ is defined modulo $\E_2$.

	For $\Phi_{14}$ we calculate for instance
	\[\Phi_{14} = -\frac{5}{4} \CC_1 - \frac{3}{4} \CC_2 - \frac{3}{2} \CC_3 - \CC_4 - \frac{1}{2} \CC_5 + \frac{1}{2} \CC_7. \]

	Now we can calculate all the local intersection numbers that define
the coefficients $\delta_p$.
	Again, we will only calculate $(D_{14},D_{32})_{fin}$ as an example, 
 since all other intersection numbers are calculated in a similar way.
	\begin{align*}
	(D_{14}, D_{32})_{fin}
	= &(\BD_{14}+\Phi_{14}, \BD_{32})_{fin}\\
	= &(\BS_1 - \BS_4 + \Phi_{14},\BS_3 - \BS_2)_{fin}\\
	= &(\BS_1, \BS_3)_{fin} - (\BS_1, \BS_2)_{fin} - (\BS_4, \BS_3 )_{fin}\\
	& + (\BS_4, \BS_2)_{fin} + (\Phi_{14}, \BS_3)_{fin} - (\Phi_{14}, \BS_2)_{fin}\\
	= &0 - m_{12} \log(5) + m_{34} \log(5)\\
	&+ m_{24} \log(2) + 0 + \frac{5}{4}(\CC_1, \BD_2)_{fin} \log(2)\\
	= &0 - \log(5) - \log(5) + 5\log(2) + 0 + \frac{5}{4} \log(2)\\
	= &\frac{25}{4} \log(2) - 2 \log(5)
	\end{align*}
	
	We can collect the intersection data into the following
	
	\begin{satz}\label{sz:busch}
          Let $D_{ij}=S_{i}-S_{j}$ and $D_{kl}=S_{k}-S_{l}$ be two
          cuspidal divisors on $(E, \belyi_E)$ with no common
          component. We write $(D_{ij}, D_{kl})_{fin} = \sum_{p \text{
              prim}} \delta_p \log(p)$.  Then only the $\delta_p$ for
          $p\in \lbrace 2,5 \rbrace$ will be different from zero and
          they are given by the table
	\begin{center}
	\begin{tabular}{c|c|c|c} 
	$D_{ij}$		&	$D_{kl}$		&	$\delta_2$	&	$\delta_5$	 \\
	\hline
	$S_1-S_4$		&	$S_3-S_2$		&	$6,25$	&	$-2$	\\
	$S_1-S_3$		&	$S_4-S_2$		&	$0$	&	$-2$	\\
	$S_1-S_2$		&	$S_3-S_4$		&	$6,25$	&	$0$	\\
	$S_1-S_4$		&	$S_3-S_0$		&	$1,25$	&	$-1$	\\
	$S_1-S_3$		&	$S_4-S_0$		&	$0$	&	$-1$	\\
	$S_1-S_0$		&	$S_3-S_4$		&	$1,25$	&	$0$	\\
	$S_1-S_4$		&	$S_0-S_2$		&	$5$	&	$-1$	\\
	$S_1-S_2$		&	$S_0-S_4$		&	$5$	&	$0$	
	\end{tabular}
	\end{center}
	\end{satz}
		\begin{beweis}
                  The necessary calculations to achieve the table have
                  been explained above; a detailed calculation is in
                  \cite{busch}.
		\end{beweis}

                \begin{bem} \label{bem:lindep} The information from
                  the table in theorem \ref{sz:busch} alone will never
                  be sufficient to calculate values of scattering
                  constants, since the the value of the formula in
                  theorem \ref{satzKuehn} stays unchanged when the
                  scattering constants differ by a constant term: Take
                  for example $D_{ij}=S_{i}-S_{j}$ and
                  $D_{kl}=S_{k}-S_{l}$ then on the left hand side of
                  theorem \ref{satzKuehn} we get the sum $C_{ik}-
                  C_{il}-C_{jk}+C_{jl}$ which would equal
                  $(C_{ik}+c)-( C_{il}+c)-(C_{jk}+c)+(C_{jl}+c)$.
	\end{bem}

        \begin{bem} Most of the calculations in this chapter can be
          done algorithmically, a detailed description of such a
          refined version of Tate's algorithm can be found in
          \cite{busch2}.
\end{bem}

\section{Linear relations for $C^{\G_E}_{jk}$}

We will use the properties of scattering constants introduced in
chapter \ref{kap:ESR} to find identities of and linear dependencies
between the scattering constants of the pair $(E,\belyi_E)$ to fix the
constant mentioned in remark \ref{bem:lindep}.

\begin{prop} \label{pr:lindepSC} For the group $\G_E$ associated with
  the Belyi pair $(E,\belyi_E)$ it suffices to know the scattering
  constants $C_{14}, C_{34}$, $C_{12}$ and $C^{\G(1)}$. Then all scattering
  constants are known and we have the list: { \allowdisplaybreaks
\begin{align*}
		C_{00} &=  \frac{1}{30}\left(C^{\G(1)}-\frac{1}{\pi}\left( 14 \log(2)+6\log(5) \right)\right) \\
		C_{01}  &=  \frac{1}{30}\left(C^{\G(1)}-\frac{1}{\pi}\left( 8 \log(2)+3\log(5) \right)\right) \\
		C_{02} &=  \frac{1}{30}\left(C^{\G(1)}-\frac{1}{\pi}\left( 2 \log(2)+3\log(5) \right)\right) \\
		C_{03} &=  \frac{1}{30}\left(C^{\G(1)}-\frac{1}{\pi}\left( 8 \log(2)+3\log(5) \right)\right) \\
		C_{04} &=  \frac{1}{30}\left(C^{\G(1)}-\frac{1}{\pi}\left( 2 \log(2)+3\log(5) \right)\right) \\
		C_{11} &= \frac{1}{4} \left( \frac{1}{6}C^{\G(1)}- \frac{62}{15 \pi}\log(2) - C_{12} \right) \\
		C_{12} & \\
		C_{13} &=  \frac{1}{4} \left(\frac{1}{6}C^{\G(1)}- \frac{32}{15 \pi}\log(2) -C_{14} \right) \\
		C_{14} & \\
		C_{22} &= \frac{1}{6}C^{\G(1)}- \frac{47}{15\pi}\log(2) -4C_{12} \\
		C_{23} &= C_{14} \\
		C_{24} &= \frac{1}{6}C^{\G(1)}- \frac{17}{15 \pi}\log(2) -4C_{14}\\
		C_{33} &= \frac{1}{4} \left( \frac{1}{6}C^{\G(1)}- \frac{62}{15 \pi}\log(2) - C_{34} \right)\\
		C_{34} & \\
		C_{44} &= \frac{1}{6}C^{\G(1)}- \frac{47}{15\pi}\log(2) -4C_{34}
		\end{align*} }\noindent
	\end{prop}
		\begin{beweis}
                  First of all, remember, that the scattering matrix
                  is symmetric, i.e. for any two cusps $S$ and $S'$ we
                  have $C_{SS'}=C_{S'S}$. Hence, in the list above
                  really all scattering constants occur.
		
                  Secondly, we apply proposition \ref{pr:sumsc}. Since
                  the point $S_0$ is totally ramified the sum in
                  formula \eqref{eq:sumsc} for $S_j=S_0$ is only to be
                  taken over one single element. When we now differ
                  $S_k$ over all the cusps, we get the first 5 rows of
                  the list above. For that we have to realize that
                  the cusps have the following widths
\begin{align*}		
b_0=10, \quad b_1=8, \quad b_2=2, \quad b_3=8, \quad b_4=2.
\end{align*}
The widths are always two times the ramification index, since the
widths of all the cusps of $\G(2)$ are two.
%In the case $C_{00}$ we have to take $C_a^{\G(2)}$ in the four other cases $C_b^{\G(2)}$. 

Now we apply proposition \ref{pr:sumsc} again but in the cases not
involving $S_0$. Then we get { \allowdisplaybreaks
 			\begin{align*}
			4 C_{11}+C_{12} & = \frac{1}{6}C^{\G(1)}- \frac{62}{15 \pi}\log(2) \\
			4 C_{33}+C_{34}& =  \frac{1}{6}C^{\G(1)}- \frac{62}{15 \pi}\log(2) \\
			4 C_{13}+C_{14} & = \frac{1}{6}C^{\G(1)}- \frac{32}{15 \pi}\log(2)\\
			4 C_{13}+C_{23} & = \frac{1}{6}C^{\G(1)}- \frac{32}{15 \pi}\log(2)\\
			4 C_{41}+C_{42} & = \frac{1}{6}C^{\G(1)}- \frac{17}{15 \pi}\log(2) \\
			4 C_{23}+C_{24} & = \frac{1}{6}C^{\G(1)}- \frac{17}{15 \pi}\log(2) \\
			4 C_{43}+C_{44}& = \frac{1}{6}C^{\G(1)}- \frac{47}{15\pi}\log(2) \\ 
			4 C_{12}+C_{22}& = \frac{1}{6}C^{\G(1)}- \frac{47}{15\pi}\log(2). 
			\end{align*}}\noindent
This follows directly with the widths from above and the fact that cusps are equivalent under $\G(2)$ if and only if they have the same image under $\belyi_E$.

		From this identities the list in the proposition follows.
		\end{beweis}

Theorem \ref{sz:busch} may now be used to calculate the missing scattering constants.

\section{Proof of the main result}

We observe, that on an elliptic curve $C$ the N\'eron Tate pairing on $\Div_0(C)$ is compatible with the  N\'eron Tate pairing on the Mordell Weil group, i.e. for $S,S' \in C(\QZ)$
	\begin{align*}   \left\langle S-\mathcal O, S'-\mathcal O \right\rangle_{NT} & = \left\langle S, S' \right\rangle_{NT}.\end{align*}
	In particular, if $C(\QZ)$ is generated by one point $P$ then the values of all N\'eron Tate pairings in $C(\QZ)$ are multiples of the  N\'eron Tate height of the generator $\hei_{NT}(P)$. Thus, by bilinearity and insertion of $\mathcal O$, for $D, D'\in \Div_0(C(\QZ))$ there is a $n\in \GZ$ such that $$ \left\langle D,D' \right\rangle_{NT}= n\cdot  \hei_{NT}(P).$$ 
We can use this fact to simplify the left hand side of equation \eqref{eq:kuehn} in theorem \ref{sz:kuehn}. 

	For the elliptic curve $E$ considered in this text we have in
        the Mordell Weil group
%get for       example: Look at the divisors $S_1-S_2$ and $S_0-S_4$.  Since
$S_2=-4S_1$, $S_3=-S_1$, $S_4=4S_1$ and $S_0=\mathcal O$ (one can use e.g. the computer algebra
        system pari to calculate multiples of $S_1$). Thus we get
	\begin{align}
	\left\langle S_1-S_2, S_0 - S_4 \right\rangle_{NT} & = \left\langle S_1, S_0 \right\rangle_{NT}-\left\langle S_1,  S_4 \right\rangle_{NT}-\left\langle S_2, S_0  \right\rangle_{NT}+\left\langle S_2,  S_4 \right\rangle_{NT} \notag \\
	& =
	0-4\hei_{NT}(S_1)- 0-16 \hei_{NT}(S_1) \notag \\
	& = -20 \hei_{NT}(S_1). \label{eq:hoehen}
	\end{align}

	\begin{satz}
		\label{sz:skc01c06}
	For the group $\G_E$ associated with $(E, \belyi_E)$ we have
	\begin{align*}
		 	C_{14} & = \frac{1}{30}\left(C^{\G(1)}+\frac{1}{\pi}\left( 7\log(2)-60\hei_{NT}(S_1) \right)\right) \\ 
			C_{21} & =   \frac{1}{30}\left(C^{\G(1)}+\frac{1}{\pi}\left( 7\log(2)-15\log(5)+60\hei_{NT}(S_1) \right)\right ) \\
			C_{34} & =  \frac{1}{30}\left(C^{\G(1)}+\frac{1}{\pi}\left( 7\log(2)-15\log(5)+60\hei_{NT}(S_1) \right)\right ) 
	\end{align*}
	and $S_1=(1,4)$ is, like before, the generator of the Mordell Weil group $E(\QZ)$.
	\end{satz}
		\begin{beweis}
	We insert information from theorem \ref{sz:busch} into theorem \ref{satzKuehn}. If we use the last row from the table in theorem \ref{sz:busch}, formula \eqref{eq:kuehn} becomes  
	\begin{align}
	\left\langle S_1-S_2, S_0 - S_4 \right\rangle_{NT} =- 5\log(2)- 2\pi \left( C_{01}- C_{14}-C_{02}+C_{42} \right). 
	\label{eq:kuehnbsp}
	\end{align}
	With the list from proposition \ref{pr:lindepSC} we simplify the sum of scattering constants to
	\begin{align}
	 2\pi \left( C_{01}- C_{14}-C_{02}+C_{24} \right) =
	 \frac{\pi}{3}C^{\G(1)} -\frac{8}{3} \log(2) -10\pi  C_{14}.
	\end{align}
	As  seen in \eqref{eq:hoehen} the right hand side of equation \eqref{eq:kuehnbsp} becomes $-20\hei_{NT}(S_1)$.
	Hence, the first scattering constant is given by
		\begin{align}
	 	C_{14} = \frac{1}{30}\left(C^{\G(1)}+\frac{1}{\pi}\left( 7\log(2)-60\hei_{NT}(S_1) \right)\right).
		\end{align}
	For the last scattering constants, we repeat the calculation above with taking the second to last and the fifth row (respectively) from the table in theorem \ref{sz:busch} and the result for $C_{14}$ into account and conclude the constants to be
		\begin{align}
	 	C_{12} =   C_{34} =      \frac{1}{30}\left(C^{\G(1)}+\frac{1}{\pi}\left( 7\log(2)-15\log(5)+60\hei_{NT}(S_1) \right)\right).
		\end{align}

		\end{beweis}

	If we insert these results into the formulas of proposition \ref{pr:lindepSC} we get a description of the scattering constants for $\G_E$ in $C^{\G(1)}$, the $\log$'s of $2$ and $5$ and $\hei_{NT}(S_1)$.

	\begin{bem}
	 Good numerical approximations exist for the N\'eron Tate height. 
 By means of the
computer algebra system pari we obtain
	$$ \hei_{NT}(S_1)\approx 0.1283750629460508690621759. $$
	This leads to the following numerical approximations of scattering constants for $\G_E$:
	\begin{align*}
		C_{00} &\approx  -0.176518865559 & C_{01}  &\approx  -0.0811617456560 \\
		C_{02} &\approx -0.0370346256255 &	C_{03} &\approx  -0.0811617456560 \\
		C_{04} &\approx  -0.0370346256255  &	C_{11} &\approx  -0.168350141906\\
		C_{12} & \approx -0.0940378417200 &		C_{13} &\approx  -0.0812067899954 \\
		C_{14} & \approx -0.00134004905741 &		C_{22} &\approx -0.170651442311\\
		C_{23} &\approx  -0.00134004905741 &		C_{24} &\approx -0.100171412657\\
		C_{33} &\approx -0.168350141906 &		C_{34} &\approx -0.0940378417200 \\
		C_{44} &\approx -0.170651442311
		\end{align*}
	\end{bem}


\begin{thebibliography}{99}

	\bibitem[Be]{belyi} Bely{\u\i}, G.: Galois extensions of a maximal cyclotomic field. Izv. Akad. Nauk SSSR Ser. Mat., Vol 43, 1979, p. 267--276.

	\bibitem[Bi]{birch} Birch, B.: Noncongruence subgroups, covers and drawings. The Grothendieck theory of dessins d'enfants (Luminy, 1993), London Math. Soc. Lecture Note Ser. Vol. 200, Cambridge University Press, Cambridge, 1994, p. 25--46.

	\bibitem[Bo]{bost} Bost, J.-B.: Introduction to compact Riemann surfaces, Jacobians, and
              abelian varieties. From number theory to physics (Les Houches, 1989), Springer, Berlin, 1992, p. 64-211.

	\bibitem[Bu]{busch} Busch, V.:  Effektive Berechnungen von N\'eron-Tate-H\"ohen mittels Arakelov-Schnittzahlen.  Diplomarbeit, Universit\"at Hamburg, 2008. 
\href{http://www.math.uni-hamburg.de/home/kuehn/diplom-busch.pdf}{http://www.math.uni-hamburg.de/home/kuehn/diplom-busch.pdf}

	\bibitem[Bu2]{busch2} Busch, V.: A refined Version of the Tate algorithm. Preprint, 2009, in preparation.

	\bibitem[Cr]{cremona} Cremona, J.: Algorithms for modular elliptic curves. Cambridge University Press, Cambridge, 1997.

	\bibitem[Dr]{drinfeld} Drinfel'd, V. G.: Two theorems on modular curves. Akademija Nauk SSSR. Funkcional' nyi Analiz i ego              Prilo\v zenija 7, 1973, No 2, p. 83-84.

	\bibitem[El1]{elkies} Elkies, N.: $ABC$ implies Mordell. Internat. Math. Res. Notices, No. 7, 1991, p. 99-109.

	\bibitem[El2]{elkik} Elkik, R.: Le th\'eor\`eme de Manin-Drinfel'd. S\'eminaire sur les Pinceaux de Courbes Elliptiques (Paris, 1988), Ast\'erisque, No. 183, 1990, p. 59-67.

	\bibitem[He]{hejhal} Hejhal, D.: The Selberg Trace Formula for $\PSL (2, \RZ)$. Vol. 2, Lecture Notes in Mathematics, 1001, Springer-Verlag, Berlin, 1983. 

	\bibitem[Hr]{hriljac} Hriljac, P.: Heights and Arakelov's intersection theory. Amer. J. Math. 107, 1985, p. 23-38.

	\bibitem[Hu]{huxley} Huxley, M.: Scattering Matrices for Congruence Subgroups. Modular forms (Durham, 1983), Horwood, Chichester, 1984.

	\bibitem[Ku]{kubota} Kubota, T.: Elementary Theory of Eisenstein Series. Halsted Press, New York, 1973. 

	\bibitem[K\"u]{kuehn} K\"uhn, U.: N\'eron-Tate heights on algebraic curves and subgroups of the modular group. Manuscripta Math., Vol. 116, 2005, No. 4, p. 401-419.

	\bibitem[La]{lang} Lang, S.: Introduction to Arakelov theory. Springer-Verlag, New York, 1988. %\todo{stimmts?}

	\bibitem[Ma]{manin} Manin, Ju. I.: Parabolic points and zeta functions of modular curves. Izv. Akad. Nauk SSSR Ser. Mat. 36, 1972, p. 19--66.

	\bibitem[Po]{posingies} Posingies, A.: Belyi-Morphismen und konstante Koeffizienten von nicht-holomorphen Eisensteinreihen. Diplomarbeit, Humboldt-Universit\"at zu Berlin, 2007. 
\href{http://www.math.uni-hamburg.de/home/kuehn/diplom-posingies.pdf}{http://www.math.uni-hamburg.de/home/kuehn/diplom-posingies.pdf}

	\bibitem[Se]{serre} Serre, J.-P.: Lectures on the Mordell-Weil theorem. Aspects of Mathematics, E15, Translated from the French and edited by Martin Brown from  notes by Michel Waldschmidt, Friedr. Vieweg \& Sohn, Braunschweig, 1989.

	
	\bibitem[Si]{silverman}Silverman, J.H.: Advanced topics in the arithmetic of elliptic curves. Graduate Texts in Mathematics 151, Springer-Verlag, New York, 1994.

	\bibitem[Ve]{venkov} Venkov, A.B.: On essentially cuspidal noncongruence subgroups of $PSL(2,\GZ)$. J. Funct. Anal. 92, 1990, p. 1-7.  

\end{thebibliography}
\end{document}